\documentclass[12pt]{article}
\usepackage{amsfonts,amsbsy,amssymb,amsmath,hyperref}
\usepackage{color}
\newcounter{dummy}

\newtheorem{thm}[equation]{Theorem}
\newtheorem{prp}[equation]{Proposition}

\newtheorem{cor}[equation]{Corollary}
\newtheorem{example}[equation]{Example}

\newenvironment{proof}{\noindent{\bf Proof: }}{\hfill$\bullet$\medskip\par}

\newcounter{ilsct}

\def\[#1{\begin{equation}\label{#1}}
\def\]{\end{equation}}

\def\cases#1{\left\{\begin{array}{ll}#1\end{array}\right.}
\def\kases#1{\left\{\begin{array}{ll}#1\end{array}\right\}}

\def\even{\mathrm{even}}
\def\odd{\mathrm{odd}}

\def\N{\mathbb{N}}
\def\NN{\mathfrak{N}}
\def\Z{\mathbb{Z}}

\def\C{\mathbb{C}}
\def\H{\mathbb{H}}
\def\X{\mathbb{X}}
\def\RP{\mathbb{R}\mathrm{P}}
\def\CP{\mathbb{C}\mathrm{P}}
\def\HP{\mathbb{H}\mathrm{P}}
\def\P{\mathbb{P}}
\def\U{\mathrm{U}}
\def\SU{\mathrm{SU}}

\def\O{\mathrm{O}}
\def\SO{\mathrm{SO}}
\def\id{\mathrm{id}}
\def\Hom{\mathrm{Hom}}
\def\r){\right)}
\def\l({\left(}
\def\qmbox#1{\quad\text{#1}\quad}
\def\diag#1{\mathrm{diag}\l(#1\r)}
\def\cc#1{\overline{#1}}

\begin{document}\sloppy

\title{ON THE TOPOLOGY OF G-MANIFOLDS WITH FINITELY MANY NON-PRINCIPAL ORBITS}
\author{S. Bechtluft-Sachs, D. J. Wraith}
\date{\today}
\maketitle

\begin{abstract}
{\small We study the topology of compact manifolds with a Lie group action for which there are only finitely many non-principal orbits, and describe the possible orbit spaces which can occur. If some non-principal orbit is singular, we show that the Lie group action must have odd cohomogeneity. We pay special attention to manifolds with one and two singular orbits, and construct some infinite families of examples. To illustrate the diversity within some of these families, we also investigate homotopy types.}
\end{abstract}


\section{Introduction}
In this paper we study $G$-manifolds, that is, smooth manifolds admitting a smooth effective action from a Lie group $G$.
We will further assume that the manifolds and groups we consider are always compact.

It is well-known that for such an action there is a unique maximal orbit type: the principal orbit. Other orbits are either `exceptional', that is, non-principal but having the same dimension as a principal orbit, or singular. Singular orbits have a strictly lower dimension than a principal orbit. In this paper we will study manifolds with Lie group action for which the non-principal orbits are finite in number, and pay special attention to the case where the non-principal orbits are singular. The motivation for studying this family of manifolds comes primarily from Geometry, which we will now explain.

The most studied families of (Riemannian) manifolds are almost certainly those which display lots of symmetry. The homogeneous spaces (equipped with homogeneous metrics) are the most symmetric family of all. These are manifolds
admitting a smooth (isometric) Lie group action which is transitive. Put another way, a homogeneous space is a manifold admitting a Lie group action for which the space of orbits consists of a single point. The topology and geometry of these spaces is for the most part well-understood.

The next most symmetric family of manifolds are those which admit a smooth action from a compact Lie group for which the dimension of the space of orbits, that is, the cohomogeneity of the $G$-manifold, is one.
Such manifolds have a simple topological description. For a compact cohomogeneity one manifold, the space of orbits is either a circle or an interval. In the first case, the manifold is just a homogeneous space bundle over the circle, and all orbits are principal. In the second case, there are two non-principal orbits corresponding to the ends of the interval,
and the manifold is a union of two disc bundles for which the non-principal orbits form the zero-section.
The boundary of each disc bundle (indeed every distance sphere, given an invariant metric) is a principal orbit, and therefore a homogeneous space. The entire manifold can be described by a group diagram, involving the main group, the
principal isotropy and the two non-principal isotropy subgroups (see \cite{AA} or \cite{GZ1}).

In a previous paper \cite{BW}, the authors studied the topology and geometry of manifolds of cohomogeneity two, three and four, in the case where the $G$-action is asystatic and the singular orbits are fixed points.

The present work begins with the observation that a cohomogeneity one manifold has at most finitely many non-principal orbits, namely zero or two. Thus is seems natural to ask about manifolds of higher cohomogeneities which share the feature of having finitely many non-principal orbits. The current paper can therefore be viewed as an outgrowth of the study of cohomogeneity one manifolds.

We will show that the behaviour of cohomogeneity one manifolds is quite different from that of $G$-manifolds with finitely many non-principal orbits in higher cohomogeneity, in that the isotropy of the latter manifolds is much more tightly constrained. It is a consequence of this constraint that we are able to describe the possible orbit spaces which arise. Specifically, we prove
\begin{thm}\label{thm-main} If $M$ is a compact $G$-manifold with finitely many non-principal orbits and cohomogeneity
greater than one, then the orbit space $G\backslash M$ consists of a manifold with boundary, together with cones over each of the boundary components.
Each boundary component corresponds to a non-principal orbit. For a singular orbit the
boundary component is either a quaternionic projective space, a complex projective space, or the quotient of $\CP^{\odd}$ by the inversion involution \eqref{quatinv-def}. For an exceptional orbit, the boundary component is either an even dimensional real projective space, or a finite quotient of an odd-dimensional sphere.
\end{thm}
If one of these non-principal orbits is singular, the corresponding boundary component is even dimensional. Hence the space of orbits
must be odd-dimensional.
As this dimension is precisely the cohomogeneity, we deduce:
\begin{cor}\label{cor-nonexcoddcoh} If $M$ is a compact $G$-manifold with at least one isolated singular orbit, then the cohomogeneity of the $G$-action must be odd.
\end{cor}
We believe that this class of manifolds is both rich and interesting, and worthy of further study from both a topological and a geometric perspective. To illustrate this richness, we construct several infinite families of examples. Firstly, we study manifolds with precisely two singular orbits, as this is the situation which most closely resembles cohomogeneity one. Our results include the following Theorem:
\begin{thm} Given any two Aloff-Wallach spaces $W_{p_1,p_2}$ and $W_{q_1,q_2}$, there is an 11-dimensional $\SU(3)$-manifold $M^{11}_{p_1p_2q_1q_2}$ of cohomogeneity three with two singular orbits equal to the given Aloff-Wallach spaces.
\end{thm}
\noindent To show the diversity within this family we also prove:
\begin{thm} From within the family $M^{11}_{p_1p_2q_1q_2}$, there is an infinite
sequence of pairwise non-homotopic manifolds for which each pair of singular orbits is non-homotopic. There is also an infinite sequence of pairwise non-homotopic `doubles', that is, manifolds for which each pair of singular orbits is identical.
\end{thm}
Unlike cohomogeneity one, it is possible to have manifolds in higher cohomogeneities which have more than two, or
precisely one singular orbit. For example we prove
\begin{thm} Given compact connected Lie groups $G\supset K$, so that $N_GK/K$ has at least rank $2$, there exist $G$-manifolds in each cohomogeneity $4k-1$, $k=1,2,3,\ldots$ with precisely one non-principal orbit.
\end{thm}

This paper is laid out as follows. In section \ref{sec-structure} we establish the basic structure of compact $G$-manifolds with finitely many singular orbits, focusing on the space of orbits and on issues which arise when we try to construct examples. In section \ref{sec-2sgorb} we study manifolds with two singular orbits and construct two explicit infinite families. We show that each of these families contains infinitely many homotopy types. In section \ref{sec-1sgorb} we study the existence of unique singular orbits, and ask about the numbers of singular orbits which are possible in a given cohomogeneity.

The authors would like to thank Dmitri Alekseevsky for encouraging us to study manifolds with finitely many singular orbits, and for his subsequent help. We would also like to thank Thomas P\"uttmann for reading a preliminary draft of this paper and for his valuable comments.

\section{The structure of $G$-spaces with finitely many non-principal orbits}\label{sec-structure}
Let $M$ be a compact smooth $G$-manifold, $G$ a compact Lie group and $K\subset G$ the principal isotropy group. We assume that all non-principal orbits of the $G$-action are isolated, in the sense that all other orbits in some $G$-invariant tubular neighbourhood of the non-principal orbit are principal.

For $q \in M,$ let $Gq$ be a non principal orbit and $H=G_q$ be the isotropy at $q$. The group $H$ acts on the normal space $\nu_q(Gq,M),$ and the normal bundle of $Gq$ in $M$ is associated to the principal $H$-bundle $G\to Gq$,
$$\nu(Gq,M)\cong G\times_H\nu_q(Gq,M)\ .$$
By our assumption that the non-principal orbits are isolated, there are no non-principal orbits of $H$ on the normal sphere $S\nu_q(Gq,M)$.
In particular, all $H$-orbits on this normal sphere have the same type.
Group actions with only one orbit type on a sphere are rather restricted. We have
\begin{thm}[\cite{Br}, p 192, Thm. 6.2]\label{bredonthm}
Let $L$ be a compact Lie group of positive dimension acting locally smoothly, effectively and with one orbit type on $S^n$. Then $L$ acts transitively or freely. If $L$ acts freely, we must have $L\cong\U(1)$, $N_{\SU(2)}\U(1)$ or $\SU(2)$. (If $\dim L=0$ then $S^n\to L\backslash S^n$ is the universal covering, so $L$ must also act freely.)
\end{thm}
A non-principal isotropy group $H$ acts transitively on the normal sphere $S^r\cong S\nu_q(Gq,M)$ if and only if the cohomogeneity of the $G$-action on $M$ is one. As this situation is well-understood, let us assume that $G$ acts with cohomogeneity $\geq 2$. Under this assumption, Theorem \ref{bredonthm} has an immediate corollary for our situation.
\begin{cor}\label{cor-gpsonsps} Let $H\to\O(r+1)$ be a representation of  a compact Lie group $H$ with only one orbit type on $S^r$. Assume that $H$ does not act transitively on the sphere $S^r$. Then the kernel of this action coincides with the isotropy group $K\subset H$ (so $K\lhd H$) and we have one of the following cases:
\begin{enumerate}
\item\label{cpn} $r=2k+1$ and the action is via complex multiplication, $H/K=\U(1)$ and the quotient space $H\backslash S^r$ is $\CP^k$;
\item\label{cpn2} $r=4k+3$ and the action is via complex multiplication, $H/K=N_{\SU(2)}\U(1)$ the normaliser of the maximal torus in $\SU(2),$ and the quotient space $N\backslash S^r$ is
\[{X-def}\X^{2k+1}:=\CP^{2k+1}/\Z_2=\CP^{2k+1}/\tau \]
where $\tau$ is the involution of $\CP^{2k+1}$ given by
\[{quatinv-def} \tau([z_1:z_2:z_3:z_4\cdots])=[-\cc{z_2}:\cc{z_1}:-\cc{z_4}:\cc{z_3}\cdots]\ ;\]
\item\label{hpn} $r=4k+3$ and the action is via quaternionic multiplication, $H/K=\SU(2)$ and the quotient space is $\HP^k$;
\item\label{finite} $H/K$ is finite and acts freely on $S^r$. If $r$ is even, the quotient must be $\RP^r$.
\end{enumerate}
\end{cor}
\begin{proof}
By Theorem \ref{bredonthm} we have $H/K\cong\U(1)$, $N_{\SU(2)}\U(1)$, $\SU(2)$ or finite, and $H/K$ acts freely on $S^r$. The irreducible representations of $\U(1)$ and $\SU(2)$ of highest weight $\neq\pm 1$ are not free. An irreducible free representation of the normaliser $N_{\SU(2)}S^1$ must restrict to an irreducible representation of the quaternion group $Q_8\subset N_{\SU(2)}S^1,$ whose only such representation is its action on the quaternions. Therefore the representation is the standard representation of $S^1$ on $\C^r$, or of $N_{\SU(2)}S^1$ resp $\SU(2)$ on $\H^r=\C^{2r}$.
\par
In the last case, $L=H/K=\O(r+1)$ is finite and acts freely on $S^r$. If $\gamma\in L\setminus\{1\}$ then $\gamma$ cannot have $+1$ as an eigenvalue. If $r=2k$ is even, then $\gamma\in\O(2k+1)$ must have $-1$ as an eigenvalue. In particular, $\gamma^2$ has an eigenvalue $+1$ and therefore $\gamma^2=\id_{S^{2k}}$. It follows that $\gamma=-\hbox{id}_{S^{2k}}$.
\end{proof}

From this the quotient $\P=H\backslash S^r$ is a projective space or a discrete quotient. We have the following possibilities for $\P$:

\[{list}\begin{tabular}{c|c|c}
$\dim\P$&singular&exceptional\\\hline
$4k$&$\HP^k$, $\CP^{2k}$&$\RP^{4k}$\\\hline
$4k+2$&$\CP^{2k+1}$, $\X^{2k+1}$&$\RP^{4k+2}$\\\hline
$2k+1$&none&$L\backslash S^{2k+1}$
\end{tabular}
\]

Let $Gq_i$, $i=1\ldots,s$ be the non-principal orbits and let $H_i=G_{q_i}$. We can choose the $q_i$ in the fixed point set of the principal isotropy group $K$ so that the ineffective kernel of all slice actions is the same and coincides with the principal isotropy group $K$, i.e.
$$K=\ker\left[H_i\to\O(\nu_{q_i}(Gq_i,M))\right]$$
We have $K\lhd H_i\subset N_GK$. Define the Weyl group $W:=N_GK/K$. Let $N_i$ be $G$-invariant tubular neighbourhoods of the singular orbits and $M^0=M\setminus\bigcup_{i=1\ldots s} N_i$.
The quotient $B=G\backslash M^0$ is a manifold whose boundary is a disjoint union of quotients of spheres $\P_i$ as listed in \eqref{list}. The $G$-invariant self-diffeomorphisms of $G/K$ are precisely those maps defined by right multiplication by elements of $W$. As a consequence, $M^0$ is the total space of a $G/K$-bundle over $B$ with a global $G$-action and structure group $W$. In other words we have
$$M^0=P_W\times_W G/K$$
for some $W$-principal bundle $P_W$.

From Corollary \ref{cor-gpsonsps} we have isomorphisms
$$\alpha_i\colon L_i\rightarrow H_i/K,$$
where $L_i=\U(1),$ $N_{\SU(2)}\U(1),$ $\SU(2)$ or a finite subgroup of $\O(r_i+1)$. Therefore the $G$-invariant tubular neighbourhoods $N_i$ and their boundaries $T_i=\partial N_i$ are $G$-equivariantly diffeomorphic to $G/K$-bundles associated to the standard $L_i$-bundle $S^{r_i}\to\P_i$, where the action of $L_i$ on $G/K$ is given by $(z,gK) \mapsto gK\alpha_i(z^{-1})$ for $z \in L_i.$ We will write this as
$$T_i\cong S^{r_i}\times_{\alpha_i} G/K\qmbox{,}N_i\cong D^{{r_i}+1}\times_{\alpha_i} G/K.$$
(Equivalently, we could view $T_i$ and $N_i$ as being associated via $\alpha_i^{-1}$ to the $H_i/K$-bundle $G/K\to G/H_i=Gq_i.$) The embedding
$$\alpha_i\colon L_i\stackrel{\cong}{\longrightarrow}H_i/K\hookrightarrow W$$
induces $W$-bundle isomorphisms
$$P_W|_{\P_i}\cong S^{r_i}\times_{\alpha_i} W.$$

The bundle $P_W$ is $W$-equivariantly diffeomorphic to $$P_W\cong M^K\cap M^0$$
where $M^K$ denotes the fixed point set of $K$.
For  connected $G$-manifold $M$ with principal isotropy $K$ one has a $G$-equivariant submersion
\[{Wred-subm}M^K\times G/K\to M\qmbox{,}(x,gK)\mapsto gx\]
which gives a $G$-equivariant diffeomorphism $M^K\times_W G/K\cong M$. Note that $W$ acts freely on $M^K\times G/K$.
Thus if $M$ is a $G$-manifold with finitely many non-principal orbits and $W=N_GK$ (or a subgroup of $N_GK$ containing all the images of the $\alpha_i$ and so that $M^0$ has a $W$-structure) then we have a $G$-equivariant diffeomorphism
\[{Wred-diff}M\cong M_W\times_W G/K\ ,\]
where
\[{Wred}M_W=P_W\cup_{\partial B}\l(\bigcup_{i=1}^sD^{r_i+1}\times_{\alpha_i}W\r)\ .\]
The $W$-manifold $M_W$ has the same orbit space as the $G$-manifold $M$ and $W$ acts freely away from the non-principal orbits.

In summary, then, given compact Lie groups $K\subset G$, a compact $G$-manifold with principal isotropy $K$ of cohomogeneity $\geq 2$ with only isolated
non-principal orbits is given by the following data, up to $G$-equivariant diffeomorphism.
\begin{enumerate}
\item A $W$-principal bundle $P_W$ over a manifold $B$ with boundary
$$\partial B=\P_1\cup\ldots\cup\P_s\ ,$$
\item Groups $H_i$, $K\lhd H_i\subset G$ together with isomorphisms $\alpha_i\colon L_i\rightarrow H_i/K$ for which
so that $P_W$ over the boundary component $\P_i$ is given by
\[{redpb}P_W|_{\P_i}=S^{r_i}\times_{\alpha_i} W.\]
\end{enumerate}
We will denote the $G$-manifold obtained from this data via \eqref{Wred-diff} and \eqref{Wred} by $M(P_W,\alpha_1,\ldots,\alpha_s).$

The question of which non-principal orbits can occur for $G$-manifolds with principal isotropy group $K$ and cohomogeneity $k$ can now in principle be decided by calculations in the non-oriented cobordism group of maps into the classifying space $BW$ for $W=N_GK/K$-bundles.

Recall that two maps $f_i\colon M_i\to Y$, $i=0,1$ from $n$-manifolds $M_i$ to a topological space $Y$ are bordant if there is a map $F\colon W\to Y$ defined on a $(n+1)$-manifold $W$ with boundary $\partial W=M_0\coprod M_1$ which extends $f_0\coprod f_1$. The cobordism group $\NN_n(Y)$ of $n$-manifolds in $Y$ is the set of such bordism classes with group structure defined by disjoint union. These groups are well understood. We have isomorphisms
$$\NN_n(Y)\cong\bigoplus_{i=0}^n H_i(Y;\Z_2)\otimes\NN_{n-i}(*)$$
where $\NN_j(*)$ is the bordism group of $j$-dimensional manifolds. Also we have that $f_i\colon M_i\to Y$, $i=0,1$, are bordant if and only if all twisted Stiefel-Whitney numbers coincide, i.e. if
\[{twswn}w_I(M_0)f_0^*(y)[M_0]=w_I(M_1)f_1^*(y)[M_1]\]
for all $y\in H^k(Y;\Z_2)$ and all partitions $I$ of $n-k$. For these facts, see \cite{BT}, \cite{Stong} or \cite{Thom} for instance.

Let $\iota_i\colon\P_i\to BL_i$ be the classifying map of the standard $L_i$-bundle $S^{r_i}\to\P_i$. Then given embeddings $\alpha_i\colon L_i\hookrightarrow W$, $i=1,\ldots, s$, a manifold $M(P_W,\alpha_1,\ldots,\alpha_s)$ exists if the map
\[{bdry-gencond}f_{\alpha_1,\ldots,\alpha_s}\colon\coprod_{i}\P_i\stackrel{\coprod_{i}\iota_i}{\longrightarrow}\coprod_{i}BL_i\stackrel{\coprod_{i}B\alpha_i}{\longrightarrow}BW\] is a boundary, that is, if all its twisted Stiefel-Whitney numbers vanish.

A necessary condition for \eqref{bdry-gencond} is of course that $\coprod_{i}\P_i$ be a boundary. By \cite{M} for example, $\HP^n$ is non-oriented cobordant to $\CP^n\times\CP^n$ and $\CP^n$ is non-oriented cobordant to $\RP^n\times\RP^n$. Also $\X^{2k+1}$ is non-oriented cobordant to $\RP^2\times\HP^k$. It is well-known that the $\RP^\odd$ are boundaries and the $\RP^\even$ are not. Hence among the above, precisely $\HP^{2k+1}$, $\CP^{2k+1}$, $\X^{4k+3}$ and $\RP^{2k+1}$ are boundaries and can appear in $G$-manifolds with precisely one non-principial orbit.

We will study such extensions in the case of a single singular orbit in section \ref{sec-1sgorb}.

\section{Manifolds with precisely two singular orbits}\label{sec-2sgorb}

As noted in the Introduction, studying manifolds in higher cohomogeneities with precisely two singular orbits is of interest as it can be viewed as a direct generalisation of the most interesting case in cohomogeneity one.

As before, let $M$ denote the manifold. Removing a small tubular neighbourhood of the two singular orbits gives
a manifold with two boundary components, which is the total space of a $G/K$-bundle over some manifold $B$.
The manifold $B$ also has two boundary components $\partial B=\P_1\cup\P_2$ with $\P_1,\P_2$ as in \eqref{list}. Computing the Stiefel-Whitney numbers shows that there are the following possibilities for $\partial B$, up to interchanging the components:
\begin{align}
\label{list2sgorb-double}\P_1&=\P_2, &L_1&=L_2&\\
\label{list2sgorb-odd}\partial B&=\CP^{2k+1}\cup\X^{2k+1},k\,\odd &L_1&=S^1, L_2=N_{\SU(2)}\U(1)&\\
\label{list2sgorb-2}\partial B&=\RP^2\cup\RP^2, &L_1&=N_{\SU(2)}\U(1), L_2=\Z_2&\\
\label{list2sgorb-exc}\partial B&=L_1\backslash S^{2k+1}\cup L_2\backslash S^{2k+1}&L_1&,L_2\subset\O(2k+2)\text{ discrete}
\end{align}

The third case is the only possibility of mixing an isolated singular with an isolated exceptional orbit due to the coincidence $\X^1=\RP^2$.

The simplest case of \eqref{list2sgorb-double} is that where $B=\P\times I$ and the quotient space $G\backslash M=\Sigma\P$ is the suspension of $\P$. For a given pair $(G,K)$ the manifold $M$ is then
$$M=M(P_W(\phi), \alpha_1,\alpha_2)$$
where
$$\alpha_1,\alpha_2\colon L=L_1=L_2\longrightarrow W$$
are so that there exists an isomorphism of $W$-bundles
$$\phi\colon S_r\times_{\alpha_1} W\to S_r\times_{\alpha_2} W$$
and $P_W(\phi)$ is the mapping cylinder of $\phi$.
Equivalently, the induced maps $\P\to BL\stackrel{B\alpha_i}{\longrightarrow}BW$ must be homotopic and $P_W$ is induced from such a homotopy $\P\times I\to BW$.

We will focus on this case, in part because of its simplicity, and in part because it is the case which most closely resembles the cohomogeneity one situation. It is important to note that unlike the cohomogeneity one case, the product $\P \times I$ is not the only candidate for the manifold $B$. Given a choice of product, take any manifold without boundary of the same dimension. Now form the connected sum between the latter manifold and the product (avoiding the boundary components). The resulting manifold is clearly also a candidate for $B$.

The simplest example of this is the double
$$M(P_W(\id),\alpha,\alpha)=D^{r+1}\times_\alpha G/K\cup_\id D^{r+1}\times_\alpha G/K$$
obtained by gluing a tubular neighbourhood of a nonprincipal orbit with itself. The next example illustrates that non-doubles are possible. Let $G=\SU(3)$ and let $K$ be trivial. Then $H_1$ and $H_2$ can be any subgroups of $\SU(3)$ isomorphic to $U(1)$ or $\SU(2)$, as $N_G(K)=\SU(3)$. For $p_1,p_2,q_1,q_2\in\Z$ with $(p_1,p_2)$ and $(q_1,q_2)$ coprime, let $\alpha_1,\alpha_2\colon\U(1)\hookrightarrow\SU(3)$ be given by
\begin{align*}
\alpha_1(z)=&\diag{z^{p_1},z^{p_2},z^{-p_1-p_2}};\\
\alpha_2(z)=&\diag{z^{q_1},z^{q_2},z^{-q_1-q_2}}
\end{align*}
where $z\in\U(1)$. The resulting homogeneous spaces $\SU(3)/\alpha(\U(1))$ are the seven-dimensional Aloff-Wallach spaces $W_{p_1,p_2}$ and $W_{q_1,q_2}$. The topology of Aloff-Wallach spaces (up to homeomorphism and diffeomorphism) is well understood. See for example \cite{KS}. In particular, among the Aloff-Wallach spaces there are infinitely many homotopy types. (This follows from the fact that $H^4(W_{a,b};\Z)\cong\Z_{a^2+ab+b^2}$.) We can therefore choose $p_1,p_2,q_1,q_2$ so that $W_{p_1,p_2}$ and $W_{q_1,q_2}$ are not diffeomorphic.

\begin{thm}\label{dim11} Given any two Aloff-Wallach spaces $W_{p_1,p_2}$ and $W_{q_1,q_2}$, there is an 11-dimensional
$\SU(3)$-manifold $M^{11}_{p_1p_2q_1q_2}$ of cohomogeneity three with two singular orbits equal to the given Aloff-Wallach spaces.
\end{thm}
\begin{proof} Let $\phi\colon S^3\times_{\alpha_1}\SU(3)\to S^3\times_{\alpha_2}\SU(3)$ be any $W=\SU(3)$-bundle isomorphism. To see that such an isomorphism exists, note that $S^3\times_{\alpha_i}\SU(3)$ is an $\SU(3)$-bundle over $\CP^1=S^2$, and these bundles are classified by $\pi_2 B\SU(3)=0.$ Now set $M^{11}_{p_1p_2q_1q_2}=M(P_{\SU(3)}(\phi),\alpha_1,\alpha_2).$
\end{proof}

As non-double examples exist, this suggests investigating the conditions under which non-double examples can arise.
This is of course a very broad question. So as to give further examples, and in particular to indicate the richness of the non-double family, we study one situation in some detail. The situation in question is the case where the space of orbits is the suspension of $\CP^m$, and where the Weyl group $W=\SU(n)$.

Recall that the cohomology of $\CP^m$ is $H^*(\CP^m;\Z)=\Z[x]/x^{m+1}$, where $x$ is the first Chern class of the universal $\U(1)$-bundle $S^{2m+1}\to\CP^m$. Since maximal tori in $W$ are conjugate, any injective homomorphism $\alpha\colon\U(1)\to W$, $W=\U(n)$ or $\SU(n)$, is conjugate to
\[{def-alpha}\alpha(p)\colon\U(1)\to W\qmbox{,}z\mapsto\diag{z^{p_1},z^{p_2},\ldots,z^{p_n}}\]
for some $p=(p_1,\ldots,p_n)\in\Z^n$. The total Chern class of $\alpha(p)$ is
\[{chernalphap}c=\sum_{k=0}^n \sigma_k(p)x^k\]
where $\sigma_k(p)$ denotes the elementary symmetric polynomial of degree $k$ in $p$, (see \cite{U}, \S8). In the case $W=\SU(n)$ we have $p_1+\cdots+p_n=0$ and the first Chern class vanishes.

\begin{prp}\label{P2.3} Let $W=\U(n)$ or $\SU(n)$ and $\alpha_1,\alpha_2\colon\U(1)\to W$ be injective homomorphisms. Then the $W$-bundles $S^r\times_{\alpha_1} W$, $S^r\times_{\alpha_2} W$ over $\CP^{(r-1)/2}$ are isomorphic if and only if they have the same Chern classes.
\end{prp}
\begin{proof} In the ``stable range'' $n>m:=(r-1)/2$ this holds for general $W$-bundles (i.e. bundles not necessarily associated to the universal bundle), see \cite{T} or \cite{OSS}, p114. If $n\leq m$ then the Chern classes of $S^r\times_{\alpha(p)} W$ determine $p$ up to permutation.
\end{proof}

We now use the above analysis to construct non-double examples of $G$-manifolds with two singular orbits.
We will again take $G=\SU(n)$ and $K$ to be trivial, so that $N_G(K)=G=SU(n)$ and $W=\SU(n)$.

\begin{thm}\label{dim13} Given Aloff-Wallach spaces $W_{p_1,p_2}$ and $W_{q_1,q_2}$, there is a 13-dimensional $\SU(3)$-manifold
$M^{13}_{p_1p_2q_1q_2}$ of cohomogeneity 5, orbit space $\Sigma \CP^2$, and two singular orbits equal to the given Aloff-Wallach manifolds if and only if
$p_1^2+p_1p_2+p_2^2=q_1^2+q_1q_2+q_2^2$.
\end{thm}
\begin{proof} Since $-\sigma_2(p_1,p_2,-p_1-p_2)=p_1^2+p_1p_2+p_2^2$, the condition guarantees that the second Chern classes coincide. Since the first Chern class of a $\SU(3)$ bundle vanishes and the third is in $H^6(\CP^2;\Z)=0$, we have from Proposition \ref{P2.3} that there is a $\SU(3)$-bundle isomorphism
$$\phi\colon S^5\times_{\alpha(p_1,p_2,-p_1-p_2)}\SU(3)\to S^5\times_{\alpha(q_1,q_2,-q_1-q_2)}\SU(3).$$
Now let $M^{13}_{p_1p_2q_1q_2}:=M(P_{\SU(3)}(\phi),\alpha(p_1,p_2,-p_1-p_2),\alpha(q_1,q_2,-q_1-q_2))$.
\end{proof}

We now have two infinite families of $\SU(3)$-manifolds with precisely two singular orbits, one family in dimension 11 and the other in dimension 13. As remarked earlier, there are infinitely many homotopy types of Aloff-Wallach manifolds, so it follows that both our families contain infinitely many equivariant diffeomorphism classes. However, if we ignore equivariance, this still leaves the question of how many diffeomorphism or homeomorphism or homotopy types occur in these families.

\begin{thm}\label{hty11} From within the family $M^{11}_{p_1p_2q_1q_2}$ in Theorem \ref{dim11}, there is an infinite sequence of pairwise non-homotopy equivalent manifolds for which each pair of singular orbits is non-homotopy equivalent. There is also an infinite sequence of pairwise non-homotopy equivalent `doubles', that is, manifolds for which each pair of singular orbits is identical.
\end{thm}

Before proving this, let us state some results which will be needed in the proof. The first of these is
\begin{thm}[\cite{K}, Theorem 0.1]\label{AWthm} The Aloff-Wallach spaces $W_{p_1,p_2}$ and $W_{q_1,q_2}$ have the same homotopy type if and only if $p_1^2+p_1p_2+p_2^2=q_1^2+q_1q_2+q_2^2$ and $p_1p_2(p_1+p_2)\equiv\pm q_1q_2(q_1+q_2)\hbox{ mod } (p_1^2+p_1p_2+p_2^2).$
\end{thm}
\begin{cor}\label{AWcor} If $p_1^2+p_1p_2+p_2^2\neq q_1^2+q_1q_2+q_2^2$, then $W_{p_1,p_2}$ and $W_{q_1,q_2}$ have different homotopy types.
\end{cor}
We will also need two number-theoretic results.
\begin{thm}[\cite{S}, 3.4]\label{NT} A positive integer $n$ is representable in the form $n=a^2+ab+b^2$ with $(a,b)=1$ if and only if the following conditions hold:
\begin{enumerate}
\item if $3^t$ divides $n$ then $t\le 1$, and
\item if $r$ is prime and $r$ divides $n$, then $r\equiv 1\hbox{ mod }3$.
\end{enumerate}
\end{thm}
The second of these number-theoretic results is a classical theorem of Dirichlet about arithmetic sequences:
\begin{thm}\label{thm-dir} Given integers $a$ and $d$ with $(a,d)=1$, there exist infinitely many
natural numbers $n$ such that $a+nd$ is prime.
\end{thm}

\begin{proof}(of Theorem \ref{hty11}.) By Theorem \ref{thm-dir} we see that there is an infinite monotonically increasing sequence of primes
$r_1,r_2,r_3,..$. all of which are congruent to 1 modulo 3. From Theorem \ref{NT} we deduce that there is a sequence of integers
$a_1,b_1,a_2,b_2,..$. such that for all natural numbers $i$:
\begin{enumerate}
\item $(a_i,b_i)=1$;
\item $a_{2i-1}^2+a_{2i-1}b_{2i-1}+b_{2i-1}^2=r_{2i-1}$;
\item $a_{2i}^2+a_{2i}b_{2i}+b_{2i}^2=r_{2i-1}r_{2i}$.
\end{enumerate}
By Corollary \ref{AWcor}, we see that the resulting Aloff-Wallach spaces $W_{a_i,b_i}$ are pairwise non-homotopy equivalent.
To complete the proof, we will show that the manifolds $M_i:=M^{11}_{a_{2i-1}b_{2i-1}a_{2i}b_{2i}}$ are pairwise non-homotopy equivalent,
by showing that the fourth cohomology groups of the $M_i$ are non-isomorphic for different $i$.

We begin this analysis by observing that $M_i$ is the union of two disc bundles (specifically $D^4$-bundles over $W_{a_{2i-1}b_{2i-1}}$
respectively $W_{a_{2i}b_{2i}}$) along their common boundaries $S^2\times SU(3)$. The Mayer-Vietoris sequence for this
union includes the following portion:
$$\cdots\rightarrow H^3(W_{a_{2i-1}b_{2i-1}})\oplus H^3(W_{a_{2i}b_{2i}})\rightarrow H^3(S^2\times SU(3))\rightarrow H^4(M_i)
\rightarrow $$
$$H^4(W_{a_{2i-1}b_{2i-1}})\oplus H^4(W_{a_{2i}b_{2i}})\rightarrow H^4(S^2\times SU(3))\rightarrow\cdots$$
The cohomology of all the spaces apart from $M_i$ is well-known, (see for example \cite{KS}, page 466 for the cohomology of Aloff-Wallach
spaces). Filling in the known groups in this sequence yields the following short exact sequence:
\[{MVseq-4}0\rightarrow\Z\rightarrow H^4(M_i)\rightarrow\Z_{r_{2i-1}}\oplus\Z_{r_{2i-1}r_{2i}}\rightarrow 0.\]
Taking tensor products with $\Z_{r_l}$ is right-exact, hence
\[{MVseq-5}\rightarrow\Z_{r_l}\rightarrow H^4(M_i)\otimes_\Z\Z_{r_l}\rightarrow
\kases{
0&l<2i-1\text{ or }l>2i\\
\Z_{r_l}^2&l=2i-1\\
\Z_{r_l}&l=2i
}
\rightarrow 0\]
is exact. In particular, for $l$ odd, we have
$$\dim_{\Z_{r_{l}}}H^4(M_i)\otimes_\Z\Z_{r_l}\cases{\geq 2&l=2i-1\\\leq 1&l\neq 2i-1}\ .$$
\end{proof}

\begin{thm} Consider the family of manifolds $M^{13}_{p_1p_2q_1q_2}$ in Theorem \ref{dim13}. If $M^{13}_{abcd}$ and $M^{13}_{a'b'c'd'}$ are two members of this family for which $$r=a^2+ab+b^2=c^2+cd+d^2\neq r'=a'^2+a'b'+b'^2=c'^2+c'd'+d'^2 ,$$
then these manifolds have different homotopy types. Consequently, the family contains infinitely many homotopy types.
\end{thm}
\begin{proof}
We show that $H^4(M^{13}_{abcd})\cong\Z_{a^2+ab+b^2}$. As in the proof of Theorem \ref{hty11}, we decompose $M=M^{13}_{abcd}$ into two disc bundles, this time $D^6$-bundles over Aloff-Wallach spaces $W_{ab}$ and $W_{cd}$ and apply the Mayer-Vietoris sequence. The common boundary $X$ of these bundles is an $\SU(3)$-bundle over $\CP^2$. The relevant portion of the Mayer-Vietoris sequence is
\[{mvseq-47}\cdots\rightarrow H^3(X)\rightarrow H^4(M)\stackrel{\chi^*}{\longrightarrow} H^4(W_{ab})\oplus H^4(W_{cd})\stackrel{\pi^*+\tilde{\pi}^*}{\longrightarrow}H^4(X)\rightarrow\cdots\]

Since $X$ also is the total space of an $S^5$-bundle over $W_{ab}$ the Gysin sequence shows that the projections $\pi\colon X\rightarrow W_{ab}$, $\tilde{\pi}\colon X\rightarrow W_{cd}$ induce isomorphisms $H^4(X)\cong H^4(W_{ab})\cong H^4(W_{cd})\cong\Z_r$ and $H^3(X)\cong H^3(W_{ab})\cong 0$. Under these isomorphisms the homomorphism $\pi^*+\tilde{\pi}^*$ corresponds to the addition map $\Z_r \times \Z_r \rightarrow \Z_r$. Thus $\chi^*$ in \eqref{mvseq-47} induces an isomorphism $H^4(M)\cong\ker(\pi^*+\tilde{\pi}^*)\cong\Z_r$.
\end{proof}

\section{Manifolds with one or many non-principal orbits}\label{sec-1sgorb}
It is interesting to compare manifolds with finitely many singular orbits in cohomogeneity at least three with those of cohomogeneity one. It is well known that compact cohomogeneity one manifolds are of two basic types: those with no non-principal orbits, in which case the manifold is the total space of a bundle over $S^1$ with the principal orbit as fibre; and those with precisely two non-principal orbits. In this latter case, the space of orbits is an interval. The non-principal orbits correspond to the end-points in the orbit space. No cohomogeneity one manifold with precisely one non-principal orbit can exist, because a point is not a boundary. For similar reasons, no cohomogeneity one manifold with more than two non-principal orbits can exist. In cohomogeneity at least three, however, the situation is very different. As noted in section \ref{sec-structure}, $\HP^{2k+1}$, $\CP^{2k+1}$, $\X^{4k+3}$ and $\RP^{2k+1}$ are all boundaries. Thus in these higher cohomogeneities, unique non-principal orbits are possible.

The orbit space of a $G$-manifold $M$ with precisely one non-principal orbit must be of the form
$$M/G=B\cup_{\P} c\P$$
where $B$ is a manifold with boundary $\partial B=\P=L\backslash S^r$, $\P$ is one of manifolds listed in \eqref{list} and $c\P$ is the cone over $\P$.
As discussed at the end of section \ref{sec-structure}, the condition for a $G$-manifold with principal isotropy $K$ with only one non-principal orbit $G/H$, $K\subset H\subset G,$ to exist is that the classifying map
\[{cl-map}f_\alpha\colon\P\stackrel{\iota}{\longrightarrow}BL\stackrel{B\alpha}{\longrightarrow}BW\]
extends to a map $B\rightarrow BW$, that is, if and only if the bordism class of $f_\alpha$ vanishes. For this we need to look at the twisted Stiefel-Whitney numbers \eqref{twswn}. Note that the map $\iota\colon\P\to BL$ induces an injection in cohomology $H^q$, $q=0,1,\ldots r$, because its homotopy fibre is $S^r$.
\bigskip

We now discuss the possibilities for unique non-principal orbits in low cohomogeneities.

\subsection*{Cohomogeneity $2$}
In cohomogeneity $2$, a non-principal orbit must be exceptional.
For such a manifold we must have $\P=S^1$ and an injective homomorphism $\alpha\colon L=Z_k\rightarrow W$.
The map $\eqref{cl-map}$ is the composition
$$f_\alpha\colon S^1 \stackrel{\iota}{\longrightarrow} B\Z_k\stackrel{B\alpha}{\longrightarrow}BW,$$
and bounds if and only if it induces $0$ in $H^1(\cdot;\Z_2)=\Hom(\pi_1(\cdot) ;\Z_2)$.
Thus an injective homomorphism $\alpha\colon\Z_k\rightarrow W$ is realized (i.e. there is a manifold $M(P_W,\alpha)$) if and only if for all $\beta$ the map
$$\Z_k\stackrel{\alpha}{\longrightarrow}W\to W/W_0\stackrel{\beta}{\longrightarrow}\Z_2$$
is zero. This always holds if $k$ is odd. Actually, in this case, the $k$-fold covering $S^1\stackrel{z\mapsto z^k}{\longrightarrow}S^1$ already bounds.

\subsection*{Cohomogeneity $3$}
In cohomogeneity $3$, by \eqref{list} we can have $\P=S^2$ and $L=H/K=\U(1)$ or $\P=\X^1=\RP^2$ with $H/K=\U(1)$ or $\Z_2$. However, since $\RP^2$ does not bound a $3$-manifold, the last two cases are ruled out. The map \eqref{cl-map} bounds if and only if the map $\U(1)\to W$ induces the zero map on $H^1(\cdot,\Z_2)$. This is automatic if for instance, the fundamental group of $W$ has odd order.

\subsection*{Cohomogeneity $4$}
The non-principal orbit must be exceptional and we must have $\P=L\backslash S^3$ with $L\subset\O(4)$ finite. Since $L$ acts freely, it preserves the orientation and therefore $w_1(\P)=0$. Note that $w_2(\P)=0$ and $w_3(\P)=0,$ since these Stiefel-Whitney classes vanish for all three dimensional manifolds. It follows that the map $f_\alpha$, with $\alpha\colon L\to W$, bounds if and only if $B\alpha\colon BL\to BW$ induces $0$ in $H^3(\cdot;\Z_2)$.

\subsection*{Cohomogeneity $5$}
The non-principal orbit must be singular and we must have $\P=\HP^1=S^4$ with $L\subset\SU(2)$. Since all Stiefel-Whitney classes of $\P$ vanish, the map $f_\alpha$, with $\alpha\colon L\to W$, bounds if and only if $B\alpha\colon BSU(2)\to BW$ induces $0$ in $H^4(\cdot;\Z_2)$.

\subsection*{Cohomogeneity $6$}
The non-principal orbit must be exceptional and we must have $\P=L\backslash S^5$ with $L\subset\SO(6)$ finite, and therefore $w_1(\P)=0$. In order that the map $f_\alpha$, with $\alpha\colon L\to W$, bounds we must have that $B\alpha\colon BL\to BW$ induces $0$ in $H^5(\cdot;\Z_2)$ but this is generally not sufficient as the example $\P=\RP^5$ shows: in this case, the non-trivial Stiefel-Whitney classes of $\P$ are $w_2$ and $w_4=w_2^2$. It follows that $f_\alpha$, with $\alpha\colon\Z_2\to W$, bounds if
$B\alpha\colon B\Z_2\to BW$ induces $0$ in $H^q(\cdot;\Z_2)$ for $q=1,3,5$.

\subsection*{Cohomogeneity $7$}
The non-principal orbit must be singular and we must have $\P=\X^3$, $L=N_{\SU(2)}\U(1)$ or $\P=\CP^3$, $L=U(1)$. The cohomology of $\X^3$ is $H^*(\X^3;\Z_2)=\Z_2[x,u]/(x^3,u^2)$ with $\deg x=1$, $\deg u =4$. The non trivial Stiefel-Whitney classes are $w_1=x$, $w_2=x^2$.  It follows that the map $f_\alpha$, with $\alpha\colon N_{\SU(2)}\U(1)\to W$, bounds if and only if $B\alpha\colon BL\to BW$ induces $0$ in $H^q(\cdot;\Z_2)$, $q=4,5,6$.
\par
In the second case, $\P=\CP^3$, $L=U(1)$, all Stiefel Whitney classes of $\P$ vanish and therefore the map $f_\alpha$, with $\alpha\colon\U(1)\to W$, bounds if and only if $B\alpha\colon BL\to BW$ induces $0$ in $H^6(\cdot;\Z_2)$.

\begin{thm}\label{su2} For $k\in\N$, let $W=\U(1)$ and $c=2k$, or $W=\SU(2)$ and $c=4k-1$. Given compact Lie groups $G\supset K$ so that $W\subset N_GK/K$, there exist $G$-manifolds with cohomogeneity $c$ with precisely one non-principal orbit.
\end{thm}
\begin{proof} Because of \eqref{Wred-diff} it suffices to construct a $W$-manifold with a single non-principal orbit. In the case $W=\U(1)$, pick $p\in\N$, $p>1$, and take the join
$$M_{W}=S^{2k+1}=S^{2k-1}*S^1$$
where $\U(1)$ acts freely on $S^{2k-1}\subset\C^k$ via the standard representation, and via $\U(1)\ni z\mapsto z^p$ on $S^1$.
\par
In the case $W=\SU(2)$ put
$$M_{W'}=S^{4k+2}=S^{4k-1}*S^2$$
where $\SU(2)$ acts freely on $S^{4k-1}\subset\C^{2k}=\l(\C^2\r)^k$ via the standard representation and on $S^2=\CP^1$ with isotropy $\U(1)$.
\end{proof}
%More generally, we could have taken $\alpha\colon z\mapsto(z^k,z^l)$ with $k,l$ coprime and so that $k+l$ is even. Then the bundle $P_W|_{\CP^{2k-1}}$ as above also bounds, though not neccessarily a bundle over $\HP^{k-1}$. To see this, let $x\in H^2(\CP^{2k-1},\Z)$ be the generator. The Stiefel Whitney classes of $P_W|_{\CP^{2k-1}}$ are $w=1+(k+l)x+klx^2=1+klx^2\mod 2$ and those of the tangent bundle $w(T\CP^{2k-1})=(1+x)^{2k}\mod 2$. The ($2\times$odd)-degree Stiefel Whitney classes of $\CP^{2k-1}$ are all zero and since $w_2(P_W|_{\CP^{2k-1}})=0$ all nonzero twisted Stiefel Whitney numbers must involve an ($2\times$odd)-degree Stiefel Whitney class of $\CP^{2k-1}$.
\begin{example} For each $c\equiv 3\hbox{ mod }4$, and each $n\ge 2$, there is an $\SU(n)$-manifold of cohomogeneity $c$ with precisely one singular orbit.
\end{example}

Given that there are plentiful examples of manifolds with precisely one and two non-principal orbits, this suggests the question of which numbers of non-principal orbits are possible.

\begin{thm} For each $c\equiv 3\hbox{ mod }4$, given compact Lie groups $G \supset K$ with $\SU(2) \subset N_GK/K,$ there is a $G$-manifold of cohomogeneity $c$ with precisely $m$ singular orbits, for each $m\in\N$.
\par
For each even $c$, given compact Lie groups $G \supset K$ with $\U(1) \subset N_GK/K,$ there is a $G$-manifold of cohomogeneity $c$ with precisely $m$ exceptional orbits, for each $m\in\N$.
\par
For each $c\equiv 1\hbox{ mod }4$, given compact Lie groups $G \supset K$ with $\U(1) \subset N_GK/K,$ there is a $G$-manifold of cohomogeneity $c$ with precisely $m$ singular orbits, for each even $m\in 2\N$.
\end{thm}
\begin{proof} The first two statements follows from Theorem \ref{su2}, by taking $m$ copies of a suitable manifold with one non-principal orbit. Away
from the non-principal orbits, these manifolds are locally $G/K$-bundles. Performing fibre connected sums then gives the desired manifold. This operation can clearly be done in an equivariant manner.
\par
For the final statement statement, take $m/2$ copies of any $G$-manifold with space of orbits $\Sigma \CP^c$, as described in section 3. Now perform fibre connected sums as before. Note that although the $\CP^{even}$ are not boundaries, any even number of disjoint copies bound.
\end{proof}

\bigskip

\noindent Department of Mathematics and Statistics, National University of Ireland Maynooth, Maynooth, Co. Kildare, Ireland.
Email: stefan@maths.nuim.ie, david.wraith@nuim.ie.
\end{document}